\documentclass{amsart}
\usepackage{graphicx}
\theoremstyle{definition}

\theoremstyle{remark}

\numberwithin{equation}{section}

\begin{document}

\title[  Infinite-sample  consistent  estimations  of  parameters]{ Infinite-sample  consistent  estimations  of  parameters of the Wiener process with drift}%

\author{Levan Labadze}%
\address{Department of  Mathematics, Georgian Technical University,  Tbilisi , Georgia}%
\email{levanlabadze@yahoo.com}%

\author{Gimzer Saatashvili}%
\address{Department of  Mathematics, Gori University,  Gori , Georgia}%
\email{gimzeri75@mail.ru}%

\author{Gogi Pantsulaia}%
\address{Department of  Mathematics, Georgian Technical University,  Tbilisi , Georgia}%
\email{gogipantsulaia@yahoo.com}%

\medskip

\medskip

\thanks{}%

\dedicatory{}%

\begin{abstract} We consider  the Wiener process with drift  $$
dX_t=\mu dt +\sigma d W_t
$$ with initial  value problem $X_0=x_0$, where
 $x_0 \in R$, $ \mu \in R$ and $\sigma > 0$ are parameters.   By use values  $(z_k)_{k \in N}$ of  corresponding  trajectories  at a fixed  positive  moment  $t$,
the infinite-sample  consistent estimates  of  each  unknown parameter of  the Wiener process with drift  are constructed under an assumption  that  all another  parameters are known. Further, we propose a certain approach for estimation of  unknown parameters  $x_0,\mu,\sigma$  of  the Wiener process  with drift   by use 
the values    $(z^{(1)}_k)_{k \in N}$      and   $(z^{(2)}_k)_{k \in N}$  being  the  results  of observations  on the  $2k$-th and  $2k+1$-th  trajectories  of  the    Wiener  process with drift  at moments  $t_1$   and $t_2$ , respectively.
\end{abstract} 
\maketitle

\section{Introduction}

Following  \cite{Susanne2012}, the Wiener process with drift  is used   as a mathematical   model   described  a random motion of  a particle suspended in water which is being
bombarded by water molecules. The temperature of the water will influence the
force of the bombardment, and thus we need a parameter $\sigma$ to characterize this.
Moreover, there is a water current which drives the particle in a certain direction,
and we will assume a parameter $\mu$ to characterize the drift. To describe the
displacements of the particle, the Wiener process can be generalized to the process 
$$
dX_t=\mu dt +\sigma d W_t  \eqno (1.1)
$$
which has solution
$$X_t= x_0 + \mu t + \sigma W_t\eqno (1.2)
$$
for $X_0=x_0$(see  \cite{Susanne2012}, p.11) . It is thus normally distributed with mean $x_0+\mu t$  and variance $\sigma^2t$, as
follows from the properties of the standard Wiener process. This process has been
proposed as a simplified model for the membrane potential evolution in a neuron.

The parameters  in (1.1)-(1.2)   have the following  sense:

(i) $\mu$ represents the equilibrium or mean value supported by fundamentals (in other words, the central location)  and called  a parameter  of the drift   in Wiener model with drift ;

(ii) $\sigma$  is  a parameter of the  bombardment  force  in Wiener model with drift ;

(iii) $x_0$ is an initioal  position of the particle suspended in water  in Wiener model with drift ;

(iv)  $x_t$ is a  position of the particle suspended in water  in Wiener model with drift at moment $t>0$;

The purpose of the present paper is to introduce a new  approach which by  use  values  $(z_k)_{k \in N}$ of  corresponding  trajectories  at a fixed  positive  moment  $t$,  will allows us
to construct a consistent estimate  for  each  unknown parameter of  theWiener model with drift  under an assumption  that  all another  parameters are known. Note that analogous  problem has been considered  by  L.Labadze and G. Pantsulaia    for  Ornstein-Uhlenbeck's stochastic process ( cf.  \cite{Labadze2016}).

The rest of the present paper is the following:

In Section 2 we consider some auxiliary notions and  facts from the theory of   mathematical statistics and  probability. 

In Section 3  we present the constructions of  consistent  and  infinite-sample  consistent  estimates  for   unknown parameters   of  the Wiener model with drift.

In Section 4  we present simulations and  animations  of   the Wiener model with drift.

In Section 5  we propose a certain approach for estimation of   unknown parameters  $x_0,\mu,\sigma$  of  the Wiener process  with drift   by use 
the values    $(z^{(1)}_k)_{k \in N}$      and   $(z^{(2)}_k)_{k \in N}$  being  the  results  of observations  on the  $2k$-th and  $2k+1$-th  trajectories  of  the    Wiener  process with drift  at moments  $t_1$   and $t_2$ , respectively.

\section{Some auxiliary notions and facts}

We begin this subsection by the following definition.

\medskip

Let  $\{ \mu_{\theta} : \theta \in R\}$  be  a family  Borel probability measures in $R$.  By  $\mu_{{\theta}}^{N}$  we  denote
the $N$-power of the measure $\mu_{\theta}$  for $ \theta \in R$.

\medskip

\noindent{\bf Definition 2.1.}  A Borel
measurable function  $T_n : R^n \to R ~(n \in N)$ is called a
consistent estimator of a parameter $\theta$ (in the sense of
everywhere convergence) for the family $(\mu_{\theta}^N)_{\theta
\in R}$ if the following condition
$$
\mu_{\theta}^N (\{ (x_k)_{k \in N} :~(x_k)_{k \in N} \in R^N~\&~
\lim_{n \to \infty}T_n(x_1, \cdots, x_n)=\theta \})=1
$$
holds true  for each $\theta \in R$.

\medskip

\noindent{\bf Definition 2.2.} A Borel
measurable function  $T_n : R^n \to R ~(n \in N)$ is called a
consistent estimator of a parameter $\theta$ (in the sense of
convergence in probability) for the family
$(\mu_{\theta}^N)_{\theta \in R}$   if for every $\epsilon>0$ and
$\theta \in R$ the following condition
$$
\lim_{n \to \infty} \mu_{\theta}^N (\{ (x_k)_{k \in N} :~(x_k)_{k
\in N} \in R^N~\&~ | T_n(x_1, \cdots, x_n)-\theta|>\epsilon \})=0
$$
holds.
\medskip

\noindent{\bf Definition 2.3.}  A Borel
measurable function  $T_n : R^n \to R ~(n \in N)$ is called a
consistent estimator of a parameter $\theta$ (in the sense of
convergence in distribution ) for the family
$(\mu_{\theta}^N)_{\theta \in R}$ if for every continuous bounded
real valued function $f$ on $R$ the following condition
$$
\lim_{n \to \infty} \int_{R^N}f(T_n(x_1, \cdots, x_n))d
\mu_{\theta}^N((x_k)_{k \in N})=f(\theta)
$$
holds.
\medskip

\noindent{\bf Remark 2.1}  Following
\cite{Shiryaev80} (see, Theorem 2, p. 272), for the family
$(\mu_{\theta}^N)_{\theta \in R}$ we have:

(a) an existence of a consistent estimator of a parameter $\theta$
in the sense of everywhere convergence implies an existence of a
consistent estimator of a parameter $\theta$ in the sense of
convergence in probability;

(b) an existence of a consistent estimator of a parameter $\theta$
in the sense of convergence in probability implies an existence of
a consistent estimator of a parameter $\theta$ in the sense of
convergence in distribution.
\medskip

\noindent{\bf Definition 2.4}  Following
\cite{Ibram80}, the family $(\mu_{\theta}^N)_{\theta \in R}$ is
called  strictly separated if there exists a  family
$(Z_{\theta})_{\theta \in R}$ of Borel subsets of $R^N$ such  that

(i)~$\mu_{\theta}^N(Z_{\theta})=1$ for $\theta \in R$;

(ii)~$Z_{\theta_1} \cap Z_{\theta_2}=\emptyset$ for all different
parameters $\theta_1$ and $\theta_2$ from $R$.

(iii)~$\cup_{\theta \in R}Z_{\theta}=R^N.$

\medskip

\noindent{\bf Definition 2.5.}  Following
\cite{Ibram80}, a Borel measurable function $T : R^N \to R$ is
called  an infinite sample  consistent  estimator of a parameter $\theta$ for
the family $(\mu_{\theta}^N)_{\theta \in R}$ if the following
condition
$$
(\forall \theta)(\theta \in R \rightarrow \mu_{\theta}^N (\{
(x_k)_{k \in N}~ : ~(x_k)_{k \in N} \in R^N ~\&~ T((x_k)_{k \in
N})=\theta\})=1)
$$
holds.
\medskip

\noindent{\bf Remark 2.2.}  Note that an
existence of an infinite sample consistent estimator of a parameter $\theta$
for the family $(\mu_{\theta}^N)_{\theta \in R}$ implies that the
family $(\mu_{\theta}^N)_{\theta \in R}$ is strictly separated.
Indeed, if we  set $Z_{\theta}=\{ (x_k)_{k \in N} : (x_k)_{k \in
N} \in R^N ~\&~T((x_k)_{k \in N})=\theta\}$ for $\theta \in R$,
then all conditions in Definition 2.4  will be satisfied.

In the sequel we will need the well known fact from the
probability theory (see, for example, \cite{Shiryaev80}, p. 390).
\medskip

\noindent{\bf Lemma 2.1.}  (Kolmogorov's strong law of large numbers)
{\it Let $X_1, X_2, ...$  be a  sequence of independent identically
distributed random variables  defined on the probability space
$(\Omega, \mathcal{F},P)$. If  these random variables have a
finite expectation $m$ (i.e., $E(X_1) = E(X_2) = ... = m <
\infty$), then the following condition
$$
P(\{ \omega : \lim_{n \to \infty} n^{-1}\sum_{k=1}^nX_k(\omega)=m
\})=1
$$
holds  true. }

\section{Estimation of parameters  of  Wiener process with drift}  

\subsection{Estimation of a parameter of the  bombardment  force  $\sigma$  in Wiener model with drift }

\medskip

The purpose of the present subsection is to estimate a parameter of the  bombardment  force  $\sigma$   by water molecules acting on a particle suspended in water  under assumption that  we know  results of observations on placements of the particle at moment $t_0$, 
a parameter  of the drift $\mu$ and an initial position $x_0$.

\noindent{\bf Theorem 3.1.1}{\it~  For  $t>0$, $x_0 \in R$,  $\mu \in R$  and $\sigma>0$,   let's  $\gamma_{(t,x_0,\mu,\sigma)}$ be a Gaussian probability  measure in  $R$ with the mean  $m_t=x_0 +\mu t$ and  the variance $\sigma_t^2= \sigma^2 t$. Assuming that  parameters $t$, $x_0$ and  $\mu$  are fixed,  denote by
$\gamma_{\sigma^2}$ the measure $\gamma_{(t,x_0,\mu,\sigma)}$.
Let define the estimate $T_n : R^n  \to R$    by the following formula
$$
T_n((z_k)_{1 \le k \le n})=\frac{ \sum_{k=1}^n (z_k -x_0 -t\mu)^2 }{n t}.    \eqno  (3.1.1)
$$
Then we get
$$
\gamma_{\sigma^2}^{\infty}\{ (z_k)_{k \in N} :  (z_k)_{k \in N} \in  R^{\infty} ~\&~   \lim_{n \to \infty}T_n((z_k)_{1 \le k \le n})=\sigma^2 \}=1,\eqno  (3.1.2)
$$
provided that $T_n$ is a consistent estimator of  a parameter of the  bombardment  force  $\sigma$  in Wiener model with drift for the family  of probability measures  $(\gamma_{\sigma^2}^{\infty})_{\sigma^2>0}$. }

\begin{proof} Let's consider probability space $(\Omega, \mathcal{F},P)$, where $\Omega=R^{\infty}$, $\mathcal{F}=B(R^{\infty})$, $P=\gamma_{\sigma^2}^{\infty}$.

For $k \in N$ we consider $k$-th projection $Pr_k$   defined on  $R^{\infty}$  by
$$
Pr_k( (x_i)_{i \in N})=x_k\eqno  (3.1.3)
$$
for $(x_i)_{i \in N} \in R^{\infty}$.

It is obvious that $(Pr_k)_{k \in N}$ is  sequence of independent Gaussian random variables with   $m_t=x_0 +\mu t$ and  the variance $\sigma_t^2= \sigma^2 t$. 
It is obvious that $(\frac{(Pr_k-x_0 -t\mu )^2}{t})_{k \in N}$  is the sequence of   independent equally distributed  random variables with mean $\sigma^2$.

 By use  Kolmogorov  Strong  Law of Large numbers we get

$$\gamma_{\sigma^2}^{\infty}\{ (z_i)_{i \in N} \in  R^{\infty} ~\&~   \lim_{n \to \infty}\sum_{k=1}^n\frac{(Pr_k( (z_i)_{i \in N})-x_0 -t\mu )^2}{tn}=\sigma^2\}=1,  \eqno  (3.1.4)
$$
which implies

$$
\gamma_{\sigma^2}^{\infty}\{ (z_i)_{i \in N}  \in  R^{\infty} ~\&~   \lim_{n \to \infty}T_n((z_k)_{1 \le k \le  n})=\sigma^2 \}=1.
$$

\end{proof}

\noindent{\bf Remark 3.1.1}  By  use  Definition 2.1, Remark 2.1   and   Theorem  3.1.1  we deduce that  $T_n$ is a consistent estimator of a parameter of the  bombardment  force  $\sigma$  in the sense of
convergence in probability for the statistical structure $(\gamma_{\sigma^2})_{\sigma^2 > 0}$  as well $T_n$ is a consistent estimator of a parameter of the  bombardment  force  $\sigma$ in the sense of
convergence in distribution for the statistical structure$(\gamma_{\sigma^2})_{\sigma^2 > 0}$. 

\medskip

\noindent{\bf Theorem 3.1.2}{\it  ~Suppose that the family of probability measures $(\gamma_{\sigma^2})_{\sigma^2 > 0}$ and the estimators $T_n:R^n \to R(n \in N)$ come  from    Theorem  4.1.1.   Then the estimators
$T^{(0)}:R^{\infty} \to  R $  and  $T^{(1)}:R^{\infty} \to  R $ defined by
$$
T^{(0)}((z_k)_{k \in N})=\underline{\lim}_{n \to \infty}T_n((z_k)_{1\le k \le n}) \eqno(3.1.5)
$$
and
$$
T^{(1)}((z_k)_{k \in N})=\overline{\lim}_{n \to \infty}T_n((z_k)_{1\le k \le n}). \eqno(3.1.6)
$$
are  infinite-sample consistent  estimators of  a parameter of the  bombardment  force  $\sigma$  in Wiener model with drift for the family  of probability measures  $(\gamma_{\sigma^2}^{\infty})_{\sigma^2>0}$.}

\begin{proof} Note that  we have
\begin{align*}
&\gamma_{\sigma^2}^{\infty}\{ (z_k)_{k \in N}  \in  R^{\infty} ~\&~  T^{(0)}((z_k)_{k \in N})=\sigma^2\}\\
&=\gamma_{\sigma^2}^{\infty}\{ (z_k)_{k \in N}  \in  R^{\infty} ~\&~ \underline{\lim}_{n \to \infty}T_n((z_k)_{1 \le k \le n})=\sigma^2 \}\\
& \ge \gamma_{\sigma^2}^{\infty}\{ (z_k)_{k \in N}  \in  R^{\infty} ~\&~   \lim_{n \to \infty}T_n((z_k)_{1 \le k \le n})=\sigma^2 \}=1,\\
\end{align*}
which means  that  $T^{(0)}$ is an infinite-sample consistent  estimator of  a parameter of the  bombardment  force  $\sigma$  in Wiener model with drift for the family  of probability measures  $(\gamma_{\sigma^2}^{\infty})_{\sigma^2>0}$.

Similarly, we have

\begin{align*}
&\gamma_{\sigma^2}^{\infty}\{ (z_k)_{k \in N} \in  R^{\infty} ~\&~  T^{(1)}((z_k)_{k \in N})=\sigma^2 \}\\
&=\gamma_{\sigma^2}^{\infty}\{ (z_k)_{k \in N} \in  R^{\infty} ~\&~ \overline{\lim}_{n \to \infty}T_n((z_k)_{1 \le k \le n})=\sigma^2 \}\\
& \ge \gamma_{\sigma^2}^{\infty}\{ (z_k)_{k \in N}  \in  R^{\infty} ~\&~   \lim_{n \to \infty}T_n((z_k)_{(z_k)_{1 \le k \le n}})=\sigma^2 \}=1,
\end{align*}
which means  that  $T^{(1)}$ is an infinite-sample consistent  estimator of  a parameter of the  bombardment  force  $\sigma$  in Wiener model with drift for the family  of probability measures  $(\gamma_{\sigma^2}^{\infty})_{\sigma^2>0}$.

\end{proof}

\noindent{\bf Remark 3.1.2}  By use  Remark 2.2  we deduce  that an
existence of  infinite sample consistent estimators  $T^{(0)}$  and  $T^{(1)}$ of a parameter of the  bombardment  force  $\sigma$  in Wiener model with drift for the family  of probability measures  $(\gamma_{\sigma^2}^{\infty})_{\sigma^2>0}$ (cf. Theorem 3.1.2) implies that the
family $(\gamma_{\sigma^2}^{\infty})_{\sigma^2>0}$ is strictly separated.

\subsection{Estimation of a parameter  of the drift $\mu$ in Wiener model with drift }

\medskip

The purpose of the present subsection is to estimate  a parameter  of the drift $\mu$ in Wiener model with drift   under assumption  that   we know  results of observations on placements of the particle at moment $t_0$, a parameter of the  bombardment  force  $\sigma$   by water molecules acting on a particle suspended in water and an initial position $x_0$.

\noindent{\bf Theorem 3.2.1}{\it~  For  $t>0$, $x_0 \in R$,  $\mu \in R$  and $\sigma>0$,   let's  $\gamma_{(t,x_0,\mu,\sigma)}$ be a Gaussian probability  measure in  $R$ with the mean  $m_t=x_0 +\mu t$ and  the variance $\sigma_t^2= \sigma^2 t$. Assuming that  parameters $t$, $x_0$ and  $\sigma$  are fixed, for $\mu \in R$  denote by
$\gamma_{\mu}$ the measure $\gamma_{(t,x_0,\mu,\sigma)}$.
Let define the estimate $T^*_n : R^n  \to R$    by the following formula
$$
T^*_n((z_k)_{1 \le k \le n})=\frac{ \sum_{k=1}^n (z_k -x_0)}{nt}.    \eqno  (3.2.1)
$$
Then we get
$$
\gamma_{\mu}^{\infty}\{ (z_k)_{k \in N} :  (z_k)_{k \in N} \in  R^{\infty} ~\&~   \lim_{n \to \infty}T^*_n((z_k)_{1 \le k \le n})=\mu \}=1,\eqno  (3.2.2)
$$
for $\mu \in R$  provided that $T_n$ is a consistent estimator of   a parameter  of the drift $\mu$ in Wiener model with drift for the family  of probability measures  $(\gamma_{\mu}^{\infty})_{\mu \in R }$. }

\begin{proof} Let's consider probability space $(\Omega, \mathcal{F},P)$, where $\Omega=R^{\infty}$, $\mathcal{F}=B(R^{\infty})$, $P=\gamma_{\sigma^2}^{\infty}$.

For $k \in N$ we consider $k$-th projection $Pr_k$   defined on  $R^{\infty}$  by
$$
Pr_k( (x_i)_{i \in N})=x_k\eqno  (3.2.3)
$$
for $(x_i)_{i \in N} \in R^{\infty}$.

It is obvious that $(Pr_k)_{k \in N}$ is  sequence of independent Gaussian random variables with   $m_t=x_0 +\mu t$ and  the variance $\sigma_t^2= \sigma^2 t$. 
It is obvious that $(\frac{Pr_k-x_0} {t})_{k \in N}$  is the sequence of   independent equally distributed  random variables with mean $\mu$.

 By use  Kolmogorov  Strong  Law of Large numbers we get

$$\gamma_{\mu}^{\infty}\{ (z_i)_{i \in N} \in  R^{\infty} ~\&~   \lim_{n \to \infty}\sum_{k=1}^n\frac{Pr_k-x_0} {tn}=\mu\}=1,  \eqno  (3.2.4)
$$
which implies

$$
\gamma_{\mu}^{\infty}\{ (z_i)_{i \in N}  \in  R^{\infty} ~\&~   \lim_{n \to \infty}T^*_n((z_k)_{1 \le k \le  n})=\mu \}=1.
$$

\end{proof}

\noindent{\bf Remark 3.2.1}  By  use  Definition 2.1, Remark 2.1   and   Theorem  3.2.1  we deduce that  $T^*_n$ is a consistent estimator of  a parameter  of the drift $\mu$ in Wiener model with drift   in the sense of
convergence in probability for the statistical structure $(\gamma_{\mu})_{\mu \in R}$  as well $T^*_n$ is a consistent estimator of  a parameter  of the drift $\mu$ in Wiener model with drift  in the sense of
convergence in distribution for the statistical structure $(\gamma_{\mu})_{\mu \in R}$.  

\medskip

\noindent{\bf Theorem 3.2.2}{\it  ~Suppose that the family of probability measures $(\gamma_{\mu})_{\mu \in R}$  and the estimators $T^*_n:R^n \to R(n \in N)$ come  from    Theorem  4.2.1.   Then the estimators
$T_*^{(0)}:R^{\infty} \to  R $  and  $T_*^{(1)}:R^{\infty} \to  R $ defined by
$$
T^{(0)}((z_k)_{k \in N})=\underline{\lim}_{n \to \infty}T^*_n((z_k)_{1\le k \le n}) \eqno(3.2.5)
$$
and
$$
T_*^{(1)}((z_k)_{k \in N})=\overline{\lim}_{n \to \infty}T^*_n((z_k)_{1\le k \le n}). \eqno(3.2.6)
$$
are  infinite-sample consistent  estimators of   a parameter  of the drift $\mu$ in Wiener model with drift  for the family  of probability measures  $(\gamma_{\mu})_{\mu \in R}$.}

\begin{proof} Note that  we have
\begin{align*}
&\gamma_{\mu}^{\infty}\{ (z_k)_{k \in N}  \in  R^{\infty} ~\&~  T_*^{(0)}((z_k)_{k \in N})=\mu\}\\
&=\gamma_{\mu}^{\infty}\{ (z_k)_{k \in N}  \in  R^{\infty} ~\&~ \underline{\lim}_{n \to \infty}T^*_n((z_k)_{1 \le k \le n})=\mu \}\\
& \ge \gamma_{\mu}^{\infty}\{ (z_k)_{k \in N}  \in  R^{\infty} ~\&~   \lim_{n \to \infty}T^*_n((z_k)_{1 \le k \le n})=\mu \}=1,\\
\end{align*}
which means  that  $T_*^{(0)}$ is an infinite-sample consistent  estimator of  a parameter  of the drift $\mu$ in Wiener model with drift  for the family  of probability measures  $(\gamma_{\mu}^{\infty})_{\mu \in R}$.

Similarly, we have

\begin{align*}
&\gamma_{\mu}^{\infty}\{ (z_k)_{k \in N} \in  R^{\infty} ~\&~  T_*^{(1)}((z_k)_{k \in N})=\mu\}\\
&=\gamma_{\mu}^{\infty}\{ (z_k)_{k \in N} \in  R^{\infty} ~\&~ \overline{\lim}_{n \to \infty}T^*_n((z_k)_{1 \le k \le n})=\mu\}\\
& \ge \gamma_{\mu}^{\infty}\{ (z_k)_{k \in N}  \in  R^{\infty} ~\&~   \lim_{n \to \infty}T^*_n((z_k)_{(z_k)_{1 \le k \le n}})=\mu \}=1,
\end{align*}
which means  that  $T_*^{(1)}$ is an infinite-sample consistent  estimator of a parameter  of the drift $\mu$ in Wiener model with drift  for the family  of probability measures  $(\gamma_{\mu}^{\infty})_{\mu \in R}$.

\end{proof}

\noindent{\bf Remark 3.2.2}  By use  Remark 2.2  we deduce  that an
existence of  infinite sample consistent estimators  $T_*^{(0)}$  and  $T_*^{(1)}$ of a parameter  of the drift $\mu$ in Wiener model with drift for the family  of probability measures  $(\gamma_{\mu}^{\infty})_{\mu \in R}$ (cf. Theorem 3.2.1) implies that the
family $(\gamma_{\mu}^{\infty})_{\mu \in R}$ is strictly separated.

\subsection{Estimation of an initioal  position of the particle suspended in water  in Wiener model with drift }

\medskip

The purpose of the present subsection is to estimate an initioal  position of the particle suspended in water under assumption that  we know results of observations on placements of the particle at moment $t_0$,  a  parameters  of the drift $\mu$ and  a parameter  of   the  bombardment  force  $\sigma$   by water molecules acting on a particle. 

\noindent{\bf Theorem 3.3.1}{\it~  For  $t>0$, $x_0 \in R$,  $\mu \in R$  and $\sigma>0$,   let's  $\gamma_{(t,x_0,\mu,\sigma)}$ be a Gaussian probability  measure in  $R$ with the mean  $m_t=x_0 +\mu t$ and  the variance $\sigma_t^2= \sigma^2 t$. Assuming that  parameters $t$, $\mu$ and  $\sigma$  are fixed,  denote by
$\gamma_{x_0}$ the measure $\gamma_{(t,x_0,\mu,\sigma)}$.
Let define the estimate $T^{**}_n : R^n  \to R$    by the following formula
$$
T^{**}_n((z_k)_{1 \le k \le n})=\frac{ \sum_{k=1}^n (z_k -t\mu)^2 }{n}.    \eqno  (3.3.1)
$$
Then we get
$$
\gamma_{x_0}^{\infty}\{ (z_k)_{k \in N} :  (z_k)_{k \in N} \in  R^{\infty} ~\&~   \lim_{n \to \infty}T^{**}_n((z_k)_{1 \le k \le n})=x_0 \}=1,\eqno  (3.3.2)
$$
for $x_0 \in R$  provided that $T^{**}_n$ is a consistent estimator of  an initioal  position of the particle suspended in water  in Wiener model with drift for the family  of probability measures  $(\gamma_{x_0}^{\infty})_{x_0 \in R}$. }

\begin{proof} Let's consider probability space $(\Omega, \mathcal{F},P)$, where $\Omega=R^{\infty}$, $\mathcal{F}=B(R^{\infty})$, $P=\gamma_{x_0}^{\infty}$.

For $k \in N$ we consider $k$-th projection $Pr_k$   defined on  $R^{\infty}$  by
$$
Pr_k( (x_i)_{i \in N})=x_k\eqno  (3.3.3)
$$
for $(x_i)_{i \in N} \in R^{\infty}$.

It is obvious that $(Pr_k)_{k \in N}$ is  sequence of independent Gaussian random variables with   $m_t=x_0 +\mu t$ and  the variance $\sigma_t^2= \sigma^2 t$. 
It is obvious that $(Pr_k-t\mu )_{k \in N}$  is the sequence of   independent equally distributed  random variables with mean $x_0$.

 By use  Kolmogorov  Strong  Law of Large numbers we get

$$\gamma_{x_0}^{\infty}\{ (z_i)_{i \in N} \in  R^{\infty} ~\&~   \lim_{n \to \infty}\sum_{k=1}^n\frac{(Pr_k( (z_i)_{i \in N})-t\mu )^2}{n}=x_0\}=1,  \eqno  (3.3.4)
$$
which implies

$$
\gamma_{x_0}^{\infty}\{ (z_i)_{i \in N}  \in  R^{\infty} ~\&~   \lim_{n \to \infty}T^{**}((z_k)_{1 \le k \le  n})=\sigma^2 \}=1.
$$

\end{proof}

\noindent{\bf Remark 3.3.1}  By  use  Definition 3.1, Remark 2.1   and   Theorem  3.3.1  we deduce that  $T^{**}$ is a consistent estimator of  an initioal  position $x_0$ of the particle suspended in water  in Wiener model with drift in the sense of
convergence in probability for the statistical structure $(\gamma_{x_0})_{x_0 \in R}$  as well $T^{**}$ is a consistent estimator ofan initioal  position $x_0$ in the sense of
convergence in distribution for the statistical structure$(\gamma_{x_0})_{x_0 \in R}$. 

\medskip

\noindent{\bf Theorem 3.3.2}{\it  ~Suppose that the family of probability measures  $(\gamma_{x_0})_{x_0 \in R}$  and the estimators $T^{**}:R^n \to R(n \in N)$ come  from    Theorem  3.3.1.   Then the estimators
$T_{***}^{(0)}:R^{\infty} \to  R $  and  $T_{***}^{(1)}:R^{\infty} \to  R $ defined by
$$
T_{***}^{(0)}((z_k)_{k \in N})=\underline{\lim}_{n \to \infty}T^{**}((z_k)_{1\le k \le n}) \eqno(3.3.5)
$$
and
$$
T_{***}^{(1)}((z_k)_{k \in N})=\overline{\lim}_{n \to \infty}T^{**}((z_k)_{1\le k \le n}). \eqno(3.3.6)
$$
are  infinite-sample consistent  estimators of   an initioal  position $x_0$ of the particle suspended in water  in Wiener model with drift  for the family  of probability measures   $(\gamma_{x_0})_{x_0 \in R}$ .}

\begin{proof} Note that  we have
\begin{align*}
&\gamma_{x_0}^{\infty}\{ (z_k)_{k \in N}  \in  R^{\infty} ~\&~ T_{***}^{(0)}((z_k)_{k \in N})=x_0\}\\
&=\gamma_{x_0}^{\infty}\{ (z_k)_{k \in N}  \in  R^{\infty} ~\&~ \underline{\lim}_{n \to \infty}T^{**}((z_k)_{1 \le k \le n})=x_0 \}\\
& \ge \gamma_{x_0}^{\infty}\{ (z_k)_{k \in N}  \in  R^{\infty} ~\&~   \lim_{n \to \infty}T^{**}((z_k)_{1 \le k \le n})=x_0 \}=1,\\
\end{align*}
which means  that  $T_{***}^{(0)}$ is an infinite-sample consistent  estimator of an initioal  position $x_0$ of the particle suspended in water  in Wiener model with drift  for the family  of probability measures $(\gamma_{x_0})_{x_0 \in R}$ .

Similarly, we have

\begin{align*}
&\gamma_{x_0}^{\infty}\{ (z_k)_{k \in N} \in  R^{\infty} ~\&~   T_{***}^{(1)}((z_k)_{k \in N})=x_0 \}\\
&=\gamma_{x_0}^{\infty}\{ (z_k)_{k \in N} \in  R^{\infty} ~\&~ \overline{\lim}_{n \to \infty}T^{**}((z_k)_{1 \le k \le n})=x_0 \}\\
& \ge \gamma_{x_0}^{\infty}\{ (z_k)_{k \in N}  \in  R^{\infty} ~\&~   \lim_{n \to \infty}T^{**}((z_k)_{(z_k)_{1 \le k \le n}})=x_0 \}=1,
\end{align*}
which means  that  $T_{***}^{(1)}$ is an infinite-sample consistent  estimator of  an initioal  position $x_0$ of the particle suspended in water  in Wiener model with drift   for the family  of probability measures  $(\gamma_{x_0})_{x_0 \in R}$ .

\end{proof}

\noindent{\bf Remark 3.3.2}  By use  Remark 2.2  we deduce  that an
existence of  infinite sample consistent estimators  $T^{(0)}$  and  $T^{(1)}$ of  an initioal  position $x_0$ of the particle suspended in water  in Wiener model with drift for the family  of probability measures  $(\gamma_{x_0})_{x_0 \in R}$ (cf. Theorem 3.3.2) implies that the
family $(\gamma_{x_0})_{x_0 \in R}$ is strictly separated.

\section{Simulations, calculations   and  animations of the Wiener process with drift}

In this section  we present some programms  in Matlab for simulation and  animation of the Wiener process with drift.  In preparation of these programms  we have used  main approaches and technique introduced in  \cite{Stanoyevitch2005}.

The simulation  of  the Wiener process with drift   can be  obtained  as follows:
$$ x_t=x_0  +\mu t  + \sigma( d_0 t+ \sqrt{2}\sum_{n=1}^{\infty}d_n\frac{\sin \pi n t}{\pi n}) , \eqno  (4.4.1)$$
where   $(d_k)_{k \in N}$ is  realization of  independent standard  Gaussian  random  variables.

If   $(d^{(k)}_n)_{ n \in \mathcal{N}}  ~( k  \in  \mathcal{N})$  is  $k$-th  realizations of the  infinite family  of  independent standard  Gaussian  random  variables, 
then  the value of the  $k$-th  trajectory  of the Wiener process with drift at moment $t$ will be

$$z_k=x_0  +\mu t  + \sigma( d^{(k)}_0 t+ \sqrt{2}\sum_{n=1}^{\infty}d^{(k)}_n\frac{\sin \pi n t}{\pi n} )\eqno (4.4.2)$$
for each $k \in N$.

In our simulation we  use  MatLab command  {\bf random('Normal',0,1,p, q)}   which  generates  $p$  coppies $(d^{(k)}_n)_{1 \le n \le q}  ~(1 \le k \le p)$ of   realizations of  the  finite  family  of   independent standard  Gaussian  random  variables  of lenght  $q$ .

In our simulation we consider the following approximation

$$z_k=x_0  +\mu t  + \sigma( d_0 t+ \sqrt{2}\sum_{n=1}^{1000}d_n\frac{\sin \pi n t}{\pi n})\eqno (4.4.3)
$$
for  $1 \le k \le  100$.

\medskip

\begin{figure}[h]
\center{\includegraphics[width=0.8\linewidth]{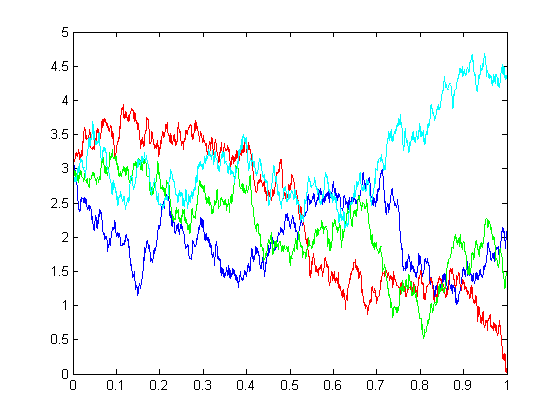}}
\caption{Four trajectories of the Wiener process with drift when $x_0=3, \mu=-1, \sigma=2.$ 
 }
\label{ris:image}
\end{figure}

The programm in MatLab giving a simulation of four trajectories of  the Wiener process with drift  with  parameters $\mu=-1$, $\sigma=3$, $x_0=3$ has the following form:

\medskip

{\bf 

$  x1 = \mbox{random}('\mbox{Normal}',0,1,4,  1001);$

$m=-1;$ 

$s=2;$ 

$x0=3;$  

$x=0 : 0.0001 : 1;$

$Y10=x0+m*x;$

$Y11=x1(1,1)*x;$

$Y20=x0+m*x;$

$Y21=x1(2,1)*x;$

$Y30=x0+m*x;$

$Y31=x1(3,1)*x;$

$Y40=x0+m*x;$

$Y41=x1(4,1)*x;$

$\mbox{for}~  k=1: 1000$

$Y11=Y11+\mbox{sqrt}(2)*x1(1,k+1)*sin(pi*k*x)/(pi*k);$

$Y21=Y21+\mbox{sqrt}(2)*x1(2,k+1)*sin(pi*k*x)/(pi*k);$

$Y31=Y31+\mbox{sqrt}(2)*x1(3,k+1)*sin(pi*k*x)/(pi*k);$

$Y41=Y41+\mbox{sqrt}(2)*x1(4,k+1)*sin(pi*k*x)/(pi*k);$

$\mbox{end}$

$X1=x;$

$Y1=Y10+s*Y11;$

$X2=x;$

$Y2=Y20+s*Y21;$

$X3=x;$

$Y3=Y30+s*Y31;$

$X4=x;$

$Y4=Y40+s*Y41;$

$\mbox{plot}(X1,Y1,'b', X2,Y2,'r-', X3,Y3,'g', X4,Y4,'c', '\mbox{LineWidth}',1)$

}

Below  we present some numerical results  obtaining by using
MatLab  and  Microsoft Excel. In our  simulation

(i) ~$n$  denotes  the number of trials;

(ii) $x_0=3$  is an initial  position of the particle suspended in water;

(iii)  $\mu=-1$  is the equilibrium or mean position supported by fundamentals;

(iv)   $\sigma=2$  is   the parameter of  bombardment  in Wiener process with drift ;

(v)   $t=0.5$ is  the  moment of the observation on the  Wiener process with drift ;

(vi)  $z_k$  is  the value of the  $k$-th  trajectory  of  Wiener process with drift    at moment $t=0.5$ (see,  Figure 1  and  Table 4.1).

\vspace{.08in} \noindent {\bf Table 4.1.  The value $z_k$  of the  the Wiener process with drift  at moment $t=0.5$   when $x_0=3, \mu=-1, \sigma=2$.     }

~\begin{center}
\begin{tabular}{| l | l | l | l | l | l | l | l | l |  l |}

\hline $ k $&  $z_k$ &      $ k $ &  $z_k$ &        $ k $ & $z_k$ &      $ k $ & $z_k$&     $ k $ & $z_k$   \\

\hline

\hline  $1 $&  $4.0991$ &      $ 21 $ &  $2.7068$ &        $ 41 $ & $4.3571$ &      $61 $ & $3.5272$&     $ 81 $ & $7.0549$  \\

\hline

\hline  $2 $&  $1.6842$ &      $ 22 $ &  $2.243$ &        $ 42 $ & $3.2793$ &      $62 $ & $2.8292$&     $ 82 $ & $2.5627$  \\

 \hline

\hline  $3 $&  $2.9422$ &      $ 23 $ &  $3.5946$ &        $ 43 $ & $3.2422$ &      $63 $ & $1.8519$&     $ 83 $ & $1.5592$  \\

 \hline

\hline  $4 $&  $4.5744$ &      $ 24 $ &  $4.6402$ &        $ 44 $ & $2.7395$ &      $64 $ & $-0.5223$&     $ 84 $ & $3.8465$  \\

 \hline

\hline  $5 $&  $2.0157$ &      $ 25 $ &  $2.2703$ &        $ 45 $ & $0.682$ &      $65 $ & $0.7602$&     $ 85 $ & $0.8471$  \\

 \hline

\hline  $6 $&  $2.8821$ &      $ 26 $ &  $2.681$ &        $ 46 $ & $3.4298$ &      $66 $ & $4.3529$&     $ 86 $ & $2.8475$  \\

 \hline

\hline  $7 $&  $4.6284$ &      $ 27 $ &  $2.9075$ &        $ 47 $ & $2.6592$ &      $67 $ & $1.1265$&     $ 87 $ & $0.5144$  \\

 \hline

\hline  $8 $&  $1.654$ &      $ 28 $ &  $3.3659$ &        $ 48 $ & $3.2093$ &      $68 $ & $1.3525$&     $ 88 $ & $2.3558$  \\

 \hline

\hline  $9 $&  $0.9561$ &      $ 29 $ &  $4.5149$ &        $ 49 $ & $0.1074$ &      $69 $ & $3.7889$&     $ 89 $ & $2.4793$  \\

 \hline

\hline  $10 $&  $-0.5407$ &      $ 30 $ &  $-1.3528$ &        $ 50 $ & $1.5462$ &      $70 $ & $3.2284$&     $ 90 $ & $2.0639$  \\

 \hline

\hline  $11 $&  $1.98941$ &      $ 31 $ &  $3.4294$ &        $ 51 $ & $1.0896$ &      $71 $ & $1.1392$&     $ 91 $ & $1.166$  \\

 \hline

\hline  $12 $&  $4.3462$ &      $ 32 $ &  $4.0825$ &        $ 52 $ & $2.1108$ &      $72 $ & $2.8833$&     $ 92 $ & $3.2738$  \\

 \hline

\hline  $13 $&  $3.455$ &      $ 33 $ &  $2.3837$ &        $ 53 $ & $2.9175$ &      $73$ & $2.3093$&     $ 93 $ & $0.4878$  \\

 \hline

\hline  $14 $&  $3.2235$ &      $ 34 $ &  $3.3037$ &        $ 54 $ & $2.9352$ &      $74 $ & $1.4444$&     $ 94 $ & $2.9504$  \\

 \hline

\hline  $15 $&  $2.1299$ &      $ 35 $ &  $4.7552$ &        $ 55 $ & $3.5531$ &      $75 $ & $-0.23$&     $ 95 $ & $2.4767$  \\

 \hline

\hline  $16 $&  $1.7554$ &      $ 36 $ &  $2.1509$ &        $ 56 $ & $2.2376$ &      $76 $ & $1.1718$&     $ 96 $ & $3.6022$  \\

 \hline

\hline  $17 $&  $3.9263$ &      $ 37 $ &  $4.3749$ &        $ 57 $ & $1.6557$ &      $77$ & $1.3588$&     $ 97 $ & $1.6005$  \\

 \hline

\hline  $18 $&  $2.4328$ &      $ 38 $ &  $2.6403$ &        $ 58 $ & $2.0076$ &      $78$ & $1.5024$&     $ 98 $ & $3.5125$  \\

 \hline

\hline  $19 $&  $2.9972$ &      $ 39 $ &  $1.4087$ &        $ 59 $ & $0.8997$ &      $79$ & $2.2622$&     $ 99 $ & $2.2026$  \\

 \hline

\hline  $20 $&  $2.59$ &      $ 40 $ &  $4.178$ &        $ 60 $ & $0.1873$ &      $80$ & $3.3423$&     $ 100 $ & $3.3646$  \\

\hline

\end{tabular}
\end{center}

\vspace{.08in} \noindent {\bf Table  4.2.  The  value of the statistic  $T_n$  for the  sample $(z_k)_{1 \le k \le n}(n=5i :1 \le i \le 20)$ from the Table 4.1.}

~\begin{center}
\begin{tabular}{|l|l|l|l|l|l|}

\hline $ n $       &    $T_n$   & $\sigma^2$ &     $n$ &   $T_n$   & $\sigma^2$\\

\hline

\hline $ 5 $       &    $3.182349256$   & $4$ &     $55$ &    $3.463530492$  & $4$\\

 \hline

\hline $ 10 $       &    $4.995431962$   & $4$ &     $60$ &   $3.472692797$ & $4$\\

\hline

\hline $ 15 $       &    $4.029170383$   & $4$ &     $65$ &   $3.628477887$   & $4$\\

\hline

\hline $ 20 $       &    $3.30673624$    & $4$ &     $70$ &   $3.621538801$   & $4$\\

\hline

\hline $ 25 $       &    $3.120603594$  & $4$ &     $75$ &   $3.662829637$   & $4$\\

\hline

\hline $ 30 $       &    $3.924002323$    & $4$ &     $80$ &   $3.554594699$  & $4$\\

\hline

\hline $ 35 $       &    $3.884200678$  & $4$ &     $85$ &   $3.96153126$   & $4$\\

\hline

\hline $ 40 $       &    $3.781846824$   & $4$ &     $90$ &   $3.836441056$   & $4$\\

 \hline

\hline $ 45 $       &    $3.715840465$   & $4$ &     $95$ &   $3.774116249$   & $4$\\

\hline

\hline $ 50 $       &    $3.665346353$   & $4$ &   $100$ &    $3.663112062$  & $4$\\

\hline

\end{tabular}
\end{center}

\medskip

\noindent{\bf Remark 4.1}  By  use  results  of  calculations   placed in the Table 4.2,  we see that the  consistent estimator $T_n$  works   successfully.

\medskip

\vspace{.08in} \noindent {\bf Table  4.3.  The  value of the statistic  $T^*_n$  for the  sample $(z_k)_{1 \le k \le n}(n=5i :1 \le i \le 20)$ from the Table 4.1.}

~\begin{center}
\begin{tabular}{|l|l|l|l|l|l|}

\hline $ n $       &    $T^*_n$   & $\mu$ &     $n$ &   $T^*_n$   & $\mu$\\

\hline

\hline $ 5 $       &    $0.12624$   & $-1$ &     $55$ &    $-0.47716$  & $-1$\\

 \hline

\hline $ 10 $       &    $-1.0209$   & $-1$ &     $60$ &   $-0.704466667$ & $-1$\\

\hline

\hline $ 15 $       &    $-0.6614$   & $-1$ &     $65$ &   $-0.851932308$   & $-1$\\

\hline

\hline $ 20 $       &    $-0.62588$    & $-1$ &     $70$ &   $-0.82396$   & $-1$\\

\hline

\hline $ 25 $       &    $-0.464312$  & $-1$ &     $75$ &   $-0.967797333$   & $-1$\\

\hline

\hline $ 30 $       &    $-0.57916$    & $-1$ &     $80$ &   $-1.0413725$  & $-1$\\

\hline

\hline $ 35 $       &    $-0.327594286$  & $-1$ &     $85$ &   $-0.959635294$   & $-1$\\

\hline

\hline $ 40 $       &    $-0.299005$   & $-1$ &     $90$ &   $-1.011635556$   & $-1$\\

 \hline

\hline $ 45 $       &    $-0.296888889$   & $-1$ &     $95$ &   $-1.056187368$   & $-1$\\

\hline

\hline $ 50 $       &    $-0.429124$   & $-1$ &   $100$ &    $-1.01773$  & $-1$\\

\hline

\end{tabular}
\end{center}

\medskip

\noindent{\bf Remark 4.2}  By  use  results  of  calculations   placed in the Table 4.2,  we see that the  consistent estimator $T^*_n$  has a tendetion will come nearer to $\mu=-1$  as soon as the number of trials 
 increases.

\medskip

The following program gives animation of the Wiener process with drift  over the time interval $[0, 1]$    when $x_0=3, \mu=-1, \sigma=2$.  

\medskip

{\bf 
$ N=100;$

$x1 = \mbox{random}('\mbox{Normal}',0,1,N,  10001);$

$m=-1;$

$s=2;$

$x0=3;$

$x=0:0.0001:1;$

$Y0=x0+m*x;$

$Y1=x1(1)*x;$

$\mbox{for}~  s=1:N $

$\mbox{for}~  k=1: 1000$

$Y1=Y1+sqrt(2)*x1(s, k+1)*sin(pi*k*x)/(pi*k);$

$\mbox{end}$

$X=x;$

$Y=Y0+s*Y1;$

$\mbox{plot}~(X,Y, 'r', '\mbox{LineWidth}',1)$

$\mbox{drawnow};$

$\mbox{pause}(1);$

$\mbox{end}$

}

 An  animation   given   by this programm  applies   $N$  different trajectories which 
are  defined by $N$  coppies  of   realizations of  the  finite  family  of   independent standard  Gaussian  random  variables  of lenght  $1001$  generated by Matlab  operator  $x1 = random('Normal',0,1,N,  10001);$

\section{ Further investigations}  Suppose  that   $(z^{(1)}_k)_{k \in N}$      and   $(z^{(2)}_k)_{k \in N}$   are results  of observations  on the  $2k$-th and  $2k+1$-th  trajectories of the Wiener process with drift   at moments  $t_1$   and $t_2,$ respectively. 
Note that having  such an  information we  can   estimate  unknown parameters $x_0,\mu,\sigma$  for  Wiener process with drift .

The first step in this direction is made by the following  proposition.

\noindent{\bf Theorem 5.1}{\it~  For  $t>0$, $x_0 \in R$,  $\mu \in R$  and $\sigma>0$,   let's  $\gamma_{(t,x_0,\mu,\sigma)}$ be a Gaussian probability  measure in  $R$ with the mean  $m_t=x_0 +\mu t$ and  the variance $\sigma_t^2= \sigma^2 t$. Assuming that  the parameter  $\sigma$  is  fixed, for  $t_1>0$,   $t_2>0$,$x_0 \in R$,  $\mu \in R$  denote by
$\gamma_{(x_0,\mu,t_1) }$  and $\gamma_{(x_0,\mu,t_2) }$the measure $\gamma_{(t_1,x_0,\mu,\sigma)}$  and $\gamma_{(t_2,x_0,\mu,\sigma)}$, restectively.

We  put  $\gamma_{(x_0,\mu)}=\gamma_{(x_0,\mu,t_1) } \times \gamma_{(x_0,\mu,t_2) }$. 
Let define the estimate ${\bf  T}_n : (R^2)^n  \to R^2$    by the following formula
$$
{\bf  T}_n ((z^{(1)}_k, z^{(2)}_k)_{1 \le k \le n})=  (  \frac{ \sum_{k=1}^n( t_2z^{(1)}_k -t_1 z^{(2)}_k)} {n(t_2-t_1)},  \frac{ \sum_{k=1}^n( z^{(2)}_k - z^{(1)}_k) } {n(t_2-t_1)} ) .    \eqno  (5.1)
$$

Then we get
$$
\gamma_{(x_0,\mu)}^{\infty}\{ (z^{(1)}_k, z^{(2)}_k)_{k \in N} :  (z^{(1)}_k, z^{(2)}_k)_{k \in N} \in  (R^2)^{\infty} $$
$$~\&~   \lim_{n \to \infty}{\bf T}_n(((z^{(1)}_k, (z^{(2)}_k)_{1 \le k \le n})=(x_0,\mu) \}=1,\eqno  (5.2)
$$
for $(x_0,\mu) \in R^2$  provided that ${\bf T_n}$ is a consistent estimator of   a  parameter $ (x_0,\mu)$  for the family  of probability measures  $(\gamma_{(x_0,\mu)}^{\infty})_{(x_0,\mu) \in R^2 }$. }

\begin{proof} Let's consider probability space $(\Omega, \mathcal{F},P)$, where $\Omega=(R^2)^{\infty}$, $\mathcal{F}=B((R^2)^{\infty})$, $P=\gamma_{(x_0,\mu)}^{\infty}$.

For $k \in N$ we consider $k$-th projection ${\bf Pr}_k= (Pr^{(1)}_k, Pr^{(2)}_k) $   defined on  $(R^2)^{\infty}$  by
$$
{\bf Pr}_k( (x^{(1)}_i,x^{(2)}_i  )_{i \in N})=(x^{(1)}_k,x^{(2)}_k  )\eqno  (5.3)
$$
for $ (x^{(1)}_i,x^{(2)}_i  )_{i \in N} \in (R^2)^{\infty}$.

It is obvious that $Pr^{(1)}_k( (x^{(1)}_i,x^{(2)}_i  )_{i \in N})=x^{(1)}_k$ and  $Pr^{(2)}_k( (x^{(1)}_i,x^{(2)}_i  )_{i \in N})=x^{(2)}_k$  are  also  projection operators for $ (x^{(1)}_i,x^{(2)}_i  )_{i \in N} \in (R^2)^{\infty}$.

It is obvious that $(Pr^{(1)}_k)_{k \in N}$ is  sequence of independent one-dimensional Gaussian random variables with   expectaion equal to $x_0 +\mu t_1$ and  the variance  equal  to $\sigma^2 t_1$. Similarly,
 $(Pr^{(2)}_k)_{k \in N}$ is  sequence of independent one-dimensional Gaussian random variables with   expectaion equal to $x_0 +\mu t_2$ and  the variance  equal  to $\sigma^2 t_2$. 
 
Note that by use Kolmgorov strong law of large numbers,  we have 

$$\gamma_{(x_0,\mu,t_1)}^{\infty}(\{ (z^{(1)}_k)_{k \in N} :  (z^{(1)}_k)_{k \in N} \in  R^{\infty} 
~\&~   \lim_{n \to \infty}  \frac{ \sum_{k=1}^nz^{(1)}_k} {n}=x_0 +\mu t_1 \})=1 \eqno (5.4)
$$
and
$$ 
\gamma_{(x_0,\mu,t_2)}^{\infty}(\{ (z^{(2)}_k)_{k \in N} :  (z^{(2)}_k)_{k \in N} \in  R^{\infty} 
~\&~  \lim_{n \to \infty}   \frac{ \sum_{k=1}^n z^{(2)}_k } {n}=\mu+\mu t_2 \})=1\eqno  (5.5)
$$

We have 

$$
\gamma_{(x_0,\mu)}^{\infty}\{ (z^{(1)}_k, z^{(2)}_k)_{k \in N} :  (z^{(1)}_k, z^{(2)}_k)_{k \in N} \in  (R^2)^{\infty}
$$
$$
 ~\&~   \lim_{n \to \infty}{\bf T}_n((z^{(1)}_k, z^{(2)}_k)_{1 \le k \le n})=(x_0,\mu) \}\eqno  (5.6)
$$
$$
=\gamma_{(x_0,\mu)}^{\infty}\{ (z^{(1)}_k, z^{(2)}_k)_{k \in N} :  (z^{(1)}_k, z^{(2)}_k)_{k \in N} \in  (R^2)^{\infty} 
$$
$$~\&~   \lim_{n \to \infty}T^{(1)}_n(((z^{(1)}_k, z^{(2)}_k)_{1 \le k \le n}=x_0
~\&   \lim_{n \to \infty}T^{(2)}_n(((z^{(1)}_k, z^{(2)}_k)_{1 \le k \le n}=\mu\}
$$
$$
=\gamma_{(x_0,\mu)}^{\infty}\{ (z^{(1)}_k, z^{(2)}_k)_{k \in N} :  (z^{(1)}_k, z^{(2)}_k)_{k \in N} \in  (R^2)^{\infty} 
$$
$$~\&~   \lim_{n \to \infty}  \frac{ \sum_{k=1}^n( t_2z^{(1)}_k -t_1 z^{(2)}_k)} {n(t_2-t_1)}=x_0
~\& \lim_{n \to \infty}   \frac{ \sum_{k=1}^n( z^{(2)}_k - z^{(1)}_k) } {n(t_2-t_1)}=\mu\}
$$
$$
=\gamma_{(x_0,\mu)}^{\infty}\{ (z^{(1)}_k, (z^{(2)}_k)_{k \in N} :  (z^{(1)}_k, (z^{(2)}_k)_{k \in N} \in  (R^2)^{\infty} 
$$
$$~\&~   \lim_{n \to \infty}  \frac{ \sum_{k=1}^nz^{(1)}_k} {n}=x_0
~\& \lim_{n \to \infty}   \frac{ \sum_{k=1}^n z^{(2)}_k } {n}=\mu\}
$$

$$
=\gamma_{(x_0,\mu,t_1)}^{\infty}(\{ (z^{(1)}_k)_{k \in N} :  (z^{(1)}_k)_{k \in N} \in  (R)^{\infty} 
~\&~   \lim_{n \to \infty}  \frac{ \sum_{k=1}^nz^{(1)}_k} {n}=x_0 +\mu t_1 \}) 
$$
$$\times 
\gamma_{(x_0,\mu,t_2)}^{\infty}(\{ (z^{(2)}_k)_{k \in N} :  (z^{(2)}_k)_{k \in N} \in  (R)^{\infty} 
~\&~  \lim_{n \to \infty}   \frac{ \sum_{k=1}^n z^{(2)}_k } {n}=\mu+\mu t_2 \})
$$
(By use  (5.4)-(5.5), we have )
$$=1 \times 1=1.
$$
\end{proof}

\medskip

\noindent{\bf Theorem 5.2}{\it  ~Suppose that the family of probability measures $(\gamma_{(x_0,\mu)}^{\infty})_{(x_0,\mu) \in R^2 }$ and the estimators  ${\bf  T}_n : (R^2)^n  \to R^2$ come  from    Theorem  5.1.   Then the estimators
${\bf T}^{(0)}:R^{\infty} \to  R $  and  ${\bf  T}^{(1)}:R^{\infty} \to  R $ defined by
$$
{\bf T}^{(0)}((z^{(1)}_k ,z^{(2)}_k  )_{k \in N})=\underline{\lim}_{n \to \infty}{\bf T}_n((z^{(1)}_k ,z^{(2)}_k )_{1\le k \le n}) \eqno(5.7)
$$
and
$$
{\bf T}^{(0)}((z^{(1)}_k ,z^{(2)}_k  )_{k \in N})=\overline{\lim}_{n \to \infty}{\bf T}_n((z^{(1)}_k ,z^{(2)}_k)_{1\le k \le n}). \eqno(5.8)
$$
are  infinite-sample consistent  estimators of   a  parameter $(x_0,\mu)$   for the family  of probability measures  $(\gamma_{(x_0,\mu)}^{\infty})_{(x_0,\mu) \in R^2 }$ .}

\begin{proof} By using  (5.2), we  get 
\begin{align*}
&\gamma_{(x_0,\mu)}^{\infty}  \{z^{(1)}_k ,z^{(2)}_k  )_{k \in N}  \in  (R^2)^{\infty} ~\&~  {\bf T}^{(0)}((z^{(1)}_k ,z^{(2)}_k  )_{k \in N})=(x_0,\mu)\}\\
&=\gamma_{(x_0,\mu)}^{\infty}  \{z^{(1)}_k ,z^{(2)}_k  )_{k \in N}  \in  (R^2)^{\infty} ~\&~  \underline{\lim}_{n \to \infty}{\bf T}_n((z^{(1)}_k ,z^{(2)}_k )_{1\le k \le n})=(x_0,\mu)\}\\
& \ge\gamma_{(x_0,\mu)}^{\infty}  \{z^{(1)}_k ,z^{(2)}_k  )_{k \in N}  \in  (R^2)^{\infty} ~\&~ \lim_{n \to \infty}{\bf T}_n((z^{(1)}_k ,z^{(2)}_k )_{1\le k \le n})=(x_0,\mu)\}=1,\\
\end{align*}
which means  that  $ {\bf T}^{(0)}$ is an  infinite-sample consistent  estimators of   a  parameter $(x_0,\mu)$   for the family  of probability measures  $(\gamma_{(x_0,\mu)}^{\infty})_{(x_0,\mu) \in R^2 }$ .

Similarly, by using  (5.2), we   have

\begin{align*}
&\gamma_{(x_0,\mu)}^{\infty}  \{z^{(1)}_k ,z^{(2)}_k  )_{k \in N}  \in  (R^2)^{\infty} ~\&~  {\bf T}^{(1)}((z^{(1)}_k ,z^{(2)}_k  )_{k \in N})=(x_0,\mu)\}\\
&=\gamma_{(x_0,\mu)}^{\infty}  \{z^{(1)}_k ,z^{(2)}_k  )_{k \in N}  \in  (R^2)^{\infty} ~\&~ \overline{\lim}_{n \to \infty}{\bf T}_n((z^{(1)}_k ,z^{(2)}_k )_{1\le k \le n})=(x_0,\mu)\}\\
& \ge\gamma_{(x_0,\mu)}^{\infty}  \{z^{(1)}_k ,z^{(2)}_k  )_{k \in N}  \in  (R^2)^{\infty} ~\&~ \lim_{n \to \infty}{\bf T}_n((z^{(1)}_k ,z^{(2)}_k )_{1\le k \le n})=(x_0,\mu)\}=1,\\
\end{align*}
which means  that  $ {\bf T}^{(1)}$   is  also  an  infinite-sample consistent  estimators of   a  parameter $(x_0,\mu)$   for the family  of probability measures  $(\gamma_{(x_0,\mu)}^{\infty})_{(x_0,\mu) \in R^2 }$ .

\end{proof}

\noindent{\bf Remark  5.1}  Following  Theorems 5.1-5.2, by using   the values    $(z^{(1)}_k)_{k \in N}$      and   $(z^{(2)}_k)_{k \in N}$  being  the  results  of observations  on the  $2k$-th and  $2k+1$-th  trajectories  of  the    Wiener  process with drift  at moments  $t_1$   and $t_2$ , respectively, we can  estimate 
parameters  $x_0$ and $\mu$.  So  an estimation  of the parameter $\sigma$  is reduced  to the case   described  in Theorem 3.1.

\bibliographystyle{amsplain}

\end{document}